# SOME CONDITIONS UNDER WHICH CERTAIN TYPES OF MODULES POSSESS LOCALIZATION PROPERTY

ISMAEL AKRAY [1] , ADIL KADIR JABBAR [2] and REZA SAZEEDEH [3]

[1] Department of Mathematics, Faculty of Science, University of Soran

[2] Department of Mathematics, School of Science, Faculty of Science and Science Education, University of Sulaimani

[3] Department of Mathematics, Faculty of Science, Urmia University

[1] E-mail: ismaeelhmd@gmail.com, [2] E-mail: adilkj@gmail.com, [3] E-mail: rsazeedeh@ipm.ir

ABSTRACT

By a localization property we mean a property that preserved under localizing modules at multiplicative closed sets. The aim of this paper is to find conditions under which we can transfere a certain property of a given module to its localization at multiplicative closed sets and conversely, that means we determine those conditions which make a given property of a module as a localization property, so we give several conditions under which certain types of modules posses localization property.

Keywords : primal, supplement, hollow, lifting, coatomic, reduced.

## 1. Introduction

Let $R$ be a commutative ring with identity. Let $M$ be an $R-$module. A proper submodule $N$ of $M$ is called a small ( or a superfluous ) $R-$module [9], denoted by $N << M$ , if $N+L \neq M$ for every proper submodule $L$ of $M$ and it is called a maximal submodule of $M$ if $L$ is any submodule of $M$ such that $N \subseteq L \subseteq M$ then $N = L$ or $L = M$ and $\text{Rad}\, M =$ intersection of all maximal submodules of $M$ [10]. $M$ is called coatomic if every proper submodule of $M$ is contained in a maximal submodule and it is called a radical module if $\text{Rad}\, M = M$ and $M$ is called a reduced module if $P(M) = 0$ , where $P(M) =$ the sum of all radical submodules of $M$ [8]. We define $N : M = \{r \in R : rM \subseteq N\}$. $M$ is called a local module [4] if it has a largest proper submodule and it is shown in [11] that an $R-$module $M$ is local if and only if $\text{Rad}\, M$ is a maximal submodule of $M$ and $\text{Rad}\, M << M$ . An element $r \in R$ is called prime to $N$ if $rm \in N$ , where $m \in M$ then $m \in N$ [1]. Thus an element $r \in R$ is not prime to $N$ if $rm \in N$ , for some $m \in M - N$ . We denote the set of all elements of $R$ that are not prime to $N$ by $S(N)$, so that $S(N) = \{r \in R : rm \in N$ , for some $m \in M - N\}$ and $N$ is called a primal submodule if $S(N)$ forms an ideal of $R$ . It is known that if $R$ is a commutative ring with identity and $S$ is a multiplicative closed set in $R$ then $R_S$ is a commutative ring with identity (the localization of $R$ at $S$ )[5, p 62]. Finally, if $M$ is an $R-$module then one can easily make $M_S$ as an $R_S$ $-$module under the module operations $\frac{m}{s}+\frac{m'}{s'}=\frac{s'm+sm'}{ss'}$ and $\frac{r}{t}\frac{m}{s}=\frac{rm}{ts}$, for all $\frac{r}{t} \in R_S$ and all $\frac{m}{s},\frac{m'}{s'} \in M_S$, so that when we say $M_S$ is a module we mean $M_S$ is an $R_S$ $-$module [6, p 74].

## 2. Some Basic Preliminaries

It is known that if a submodule $N$ of an $R-$module $M$ is prime then $(N : M)$ is a prime ideal of $R$ [3] and if $N$ is a primal submodule of $M$ then $S(N)$ is a prime ideal of $R$ [2] (and hence a primal ideal of $R$ ). If $N, L$ are submodules of $M$ and $S$ is a multiplicative closed set in $R$ then one can easily show that $(L+N)_S = L_S + N_S$ and $(L \cap N)_S = L_S \cap N_S$ . Also, one can prove that $N_S$ is a submodule of $M_S$ . If $x \in M$ and $s \in S$ such that $\frac{x}{s} = 0$ then for some $q \in S$ , we have $qx = 0$ and if $r \in R, s \in S$ then $\frac{r}{s}M_S = (rM)_S$ .

Throughout this paper $R$ is a commutative ring with identity and $M$ is a unitary $R-$module.

## 3. Main Results

It is known that every strongly prime submodule is a prime submodule but the converse is not true in general ( see [7, Example 2.2] ). In this work, we generalize the result in [7, Theorem 2.3] to prime submodules.

Proposition 3.1. Let $M$ be an $R-$module and $N$ is a proper submodule of $M$ . Then $N$ is prime if and only if for every ideal $I$ of $R$ and every submodule $L$ of $M$ with $IL \subseteq N$ then $L \subseteq N$ or $I \subseteq (N : M)$ .

Proof. Let $N$ be prime and $I$ is an ideal of $R$ and $L$ is a submodule of $M$ with $IL \subseteq N$ and $L \not\subseteq N$ . Then there exists $x \in L$ and $x \notin N$. Now, for all $y \in I$ we have $yx \in IL \subseteq N$ and as $N$ is prime we get $y \in (N:M)$ and so $I \subseteq (N:M)$.

Conversely, suppose that the given hypothesis holds. Let for any $r \in R$ and $m \in M$ we have $rm \in N$. Then $\langle r \rangle$ is an ideal of $R$ and $Rm$ is a submodule of $M$. If $t \in \langle r \rangle$ and $b \in Rm$, then $t = sr$ and $b = am$ for some $s, a \in R$. So that $tb = sram = sarm \in N$, which gives that $\langle r \rangle Rm \subseteq N$. Then by the given condition we have $Rm \subseteq N$ or $\langle r \rangle \subseteq (N:M)$, the former case gives $m \in N$ and the later one gives $r \in (N:M)$. Hence $N$ is a prime submodule of $M$. ■

Corollary 3.2. Let $M$ be an $R-$module and $N$ is a proper submodule of $M$. Then $N$ is prime if and only if for any $y \in R$ and every submodule $L$ of $M$ with $yL \subseteq N$, then $L \subseteq N$ or $y \in (N:M)$.

Proof. Suppose that $N$ is prime and let $y \in R$ and $L$ is a submodule of $M$ with $yL \subseteq N$. Clearly $\langle y \rangle$ is an ideal of $R$ and $yL \subseteq N$ gives $\langle y \rangle L \subseteq N$ and as $N$ is prime by Proposition 3.1, we get $L \subseteq N$ or $\langle y \rangle \in (N:M)$, the last case gives $y \in (N:M)$.

Conversely, suppose that the given condition holds. To show that $N$ is prime. Let $I$ be any ideal of $R$ and $L$ be any submodule of $M$ such that $IL \subseteq N$. As $I \neq \phi$, take $y \in I$ and then we have $yL \subseteq N$ so by the condition we get $L \subseteq N$ or $y \in (N:M)$, this last possibility gives $I \subseteq (N:M)$ and thus by Proposition 3.1, we get $N$ is prime. ■

Proposition 3.3. Let $M$ be an $R-$module and $N$ is a proper submodule of $M$. Then we have $S(N:M) \subseteq S(N)$ ( and hence $(N:M) \subseteq S(N)$ ).

Proof. Let $r \in S(N:M)$, that means $rs \in N:M$ for some $s \notin N:M$, so that $rsM \subseteq N$ and $sM \not\subset N$, which gives that $sm \notin N$ for some $m \in M$ and then $rsm \in N$ so that $r \in S(N)$ and thus $S(N:M) \subseteq S(N)$. ■

Proposition 3.4. Let $M$ be an $R-$module. If $N$ is a primal submodule of $M$, then $(N:M)$ is a primal ideal of $R$.

Proof. Since $1 \notin (N:M)$ and $0 \in R$ such that $0.1 = 0 \in (N:M)$ so that $0 \in S(N:M)$ and thus $\phi \neq S(N:M) \subseteq R$. Let $r, s \in S(N:M)$, then $ra \in N:M$ and $sb \in (N:M)$ for some and $a, b \notin (N:M)$. Then $raM \subseteq N$, $sbM \subseteq N$, from which we get $rabM \subseteq N$ and $sbaM \subseteq N$, so that we get $(r-s)abM \subseteq N$ and that means $(r-s)ab \in (N:M)$ and as $a \notin (N:M)$ we get $(r-s)b \in S(N:M)$. If $r-s \notin S(N:M)$, then we have $r-s \notin (N:M)$, this gives $b \in (N:M)$, which is a contradiction and thus we get $r-s \in S(N:M)$. Let $c \in R$, then $rcaM \subseteq N$, so that $rca \in (N:M)$. Since $a \notin N:M$, we get $rc \in S(N:M)$ and so $S(N:M)$ is an ideal of $R$. Hence $(N:M)$ is a primal ideal of $R$. ■

Definition 3.5. Let $M$ be an $R-$module. We call a proper submodule $N$ of $M$ as a complementary submodule of $M$ if $r \in R$ is such that $rx \notin N$ for some $x \notin N$ then $ry \notin N$ for all $y \notin N$. Also, we define $N:r = \{x \in M : rx \in N\}$.

Proposition 3.6. Let $M$ be an $R-$module and $N$ is a complementary submodule of $M$ and $(N:M)$ is a primal ideal of $R$, then $N$ is a primal submodule of $M$.

Proof. To show $S(N) = S(N:M)$. Clearly $S(N:M) \subseteq S(N)$. Now, suppose that $S(N) \not\subset S(N:M)$, so that there exists $r \in S(N)$ and $r \notin S(N:M)$, then $rx \in N$ for some $x \notin N$ and $r \notin (N:M)$, so that $rm \notin N$ for some $m \in M$ and then we get $m \notin N$. As $N$ is complementary we get $rx \notin N$ which is a contradiction, so $S(N) \subseteq S(N:M)$. Hence $S(N) = S(N:M)$. Now, since $(N:M)$ is a primal ideal of $R$, so that $S(N:M)$ is an ideal of $R$ and so that $S(N)$ is an ideal of $R$ and thus $N$ is a primal submodule of $M$. ■

Combining Proposition 3.4 and Proposition 3.6, we get the following corollary.

Corollary 3.7. Let $M$ be an $R-$module and $N$ is a complementary submodule of $M$, then $N$ is a primal submodule of $M$ if and only if $(N:M)$ is a primal ideal of $R$.

Proposition 3.8. Let $R$ be a commutative ring with identity. If $\phi \neq A \subseteq R$ and $S$ is a multiplicative closed set in $R$ such that $S(A) \cap S = \phi$, then $A$ is an ideal of $R$ if and only if $A_S$ is an ideal of $R_S$.

Proof. Let $A$ be an ideal of $R$. To show $A_S$ is an ideal of $R_S$. Let $\frac{a}{s}, \frac{a'}{s'} \in A_S$ and $\frac{r}{t} \in R_S$ where $a, a' \in A, r \in R$ and $s, s', t \in S$. Then $\frac{a}{s} - \frac{a'}{s'} = \frac{s'a - sa'}{ss'} \in A_S$ and also we have $\frac{r}{t}\frac{a}{s} = \frac{ra}{ts} \in A_S$ and $\frac{a}{s}\frac{r}{t}=\frac{ar}{st}\in A_S$ , so that $A_S$ is an ideal of $R_S$ . Let $A_S$ be an ideal of $R_S$ . Let $a,b\in A$ and $r\in R$ , then for any $s\in S$ we have $\frac{a}{s},\frac{b}{s}\in A_S$ and $\frac{r}{s}\in R_S$ , then $\frac{a-b}{s}=\frac{a}{s}-\frac{b}{s}\in A_S$ , so that $\frac{a-b}{s}=\frac{c}{u}$ , for some $c\in A$ and $u\in S$ , then there exists $v\in S$ such that $vu(a-b)=vsc\in A$ . If $a-b\notin A$ then we get $vu\in S(A)$ and thus $vu\notin S$ which is a contradiction, so that $a-b\in A$ . Also we have $\frac{ra}{ss}=\frac{r}{s}\frac{a}{s}\in A_S$ , so that $\frac{ra}{ss}=\frac{a'}{s'}$ and then there exists $s''\in S$ such that $s''s'ra=s''ssa'\in A$ . If $ra\notin A$ then $s''s'\in S(A)$ and thus $s''s'\notin S$ , which is a contradiction so that $ra\in A$ . By a similar argument we get $ar\in A$ . Hence $A$ is an ideal of $R$ . ■

Corollary 3.9. Let $M$ be an $R-$module and $N$ is a submodule of $M$ . If $S$ is a multiplicative closed set in $R$ such that $S(N)\cap S=\phi$ , then $S(N)$ is an ideal of $R$ if and only if $(S(N))_S$ is an ideal of $R_S$ .

Proof. Put $S(N)=A$ in Proposition 3.8, the result follows. ■

Theorem 3.10. Let $M$ be an $R-$module and $S$ is a multiplicative system in $R$ . If $N$ is a proper submodule of $M$ such that $S(N)\cap S=\phi$ then $S(N_S)=(S(N))_S$ .

Proof. To show that $S(N_S)=(S(N))_S$ . Let $\frac{r}{s}\in S(N_S)$ , where $r\in R$ and $s\in S$ . Then $\frac{r}{s}\frac{m}{t}\in N_S$ , for some $\frac{m}{t}\notin N_S$ , which gives $\frac{rm}{st}=\frac{n}{u}$ , for some $n\in N, u\in S$ and $m\notin N$ , so that $qurm=qstn\in N$ , for some $q\in S$ . As $m\notin N$ , we get $qur\in S(N)$ and thus $\frac{r}{s}=\frac{qur}{qus}\in(S(N))_S$ . Hence we have $S(N_S)\subseteq(S(N))_S$ .

Conversely, let $\frac{r}{s}\in(S(N))_S$ , for $r\in S(N)$ and $s\in S$ . So that $rm\in N$ , for some $m\notin N$ . If $\frac{m}{s}\in N_S$ , then $\frac{m}{s}=\frac{m'}{s'}$ , for some $m'\in N$ and $s'\in S$ . Then $qs'm=qsm'\in N$ and as $m\notin N$ we get $qs'\in S(N)$ and as $S(N)\cap S=\phi$ we get $qs'\notin S$ that is a contradiction. Hence we get $\frac{m}{s}\notin N_S$ and as $\frac{r}{s}\frac{m}{s}=\frac{rm}{ss}\in N_S$ , we have $\frac{r}{s}\in S(N_S)$ , thus $(S(N))_S\subseteq S(N_S)$ . Hence $S(N_S)=(S(N))_S$ . ■

Theorem 3.11. Let $M$ be an $R-$module and $S$ is a multiplicative system in $R$ . If $N$ is a proper submodule of $M$ such that $S(N)\cap S=\phi$ then $N$ is primal if and only if $N_S$ is a primal submodule of $M_S$ .

Proof. Let $N_S=M_S$ , then if $x\in M$ we get $\frac{x}{s}\in M_S=N_S$ , for $s\in S$ , so that $\frac{x}{s}=\frac{y}{t}$ , for some $y\in N$ and $t\in S$ . Thus we have $qtx=qsy\in N$ , for some $q\in S$ . As $qt\in S$ we have $qt\notin S(N)$ and thus we get $x\in N$ . Hence $M=N$ which is a contradiction and thus $N_S\neq M_S$ , that means $N_S$ is a proper submodule of $M_S$ . Let $N$ be primal so that $S(N)$ is an ideal of $R$ and thus by Corollary 3.9, we have $(S(N))_S$ is an ideal of $R_S$ that means $S(N_S)$ is an ideal of $R_S$ and thus $N_S$ is a primal submodule of $M_S$ . By Theorem 3.10, we have $S(N_S)=(S(N))_S$ , from which we conclude that $S(N)$ is an ideal of $R$ if and only if $S(N_S)$ is an ideal of $R_S$ and so that $N$ is a primal submodule of $M$ if and only if $N_S$ is a primal submodule of $M_S$ . ■

Example 3.12. Consider the $Z-$module $Z_6$ and the submodule $N=\{0,2,4\}$ of $Z_6$ . If we take the multiplicative closed set $S=\{-1,1\}$ in $Z$ then we have $(N:-1)=\{x\in Z_6:-x\subseteq N\}=N$ and also we have $(N:1)=\{x\in Z_6:x\subseteq N\}=N$ , that means $N:s=N$ , for all $s\in S$ .

Proposition 3.13. Let $M$ be an $R-$module and $S$ is a multiplicative closed set in $R$ . If $N$ is a proper submodule of $M$ such that $N:s=N$ for all $s\in S$ then $N$ is prime if and only if $N_S$ is a prime submodule of $M_S$ .

Proof. Let $N$ be a prime submodule of $M$ . Let for $\frac{r}{s}\in R_S$ , $\frac{x}{t}\in N_S$ we have $\frac{r}{s}\frac{x}{t}\in N_S$ so that $\frac{rx}{st}=\frac{y}{u}$ for some $y\in N$ and $u\in S$ then $qurx=qsty\in N$ for some $q\in S$ and as $N$ is a prime submodule of $M$ we have $x\in N$ or $qurM\subseteq N$ from which we get $\frac{x}{t}\in N_S$ or $\frac{r}{s}M_S=\frac{qur}{qus}M_S=(qurM)_S\subseteq N_S$ and thus $N_S$ is a prime submodule of $M_S$ .

Conversely, suppose that $N_S$ is a prime submodule of $M_S$. Let $r \in R, m \in M$ be such that $rm \in N$ then for $s \in S$ we have $\frac{r}{s}\frac{m}{s} \in N_S$ and as $N_S$ is prime we get $\frac{m}{s} \in N_S$ or $\frac{r}{s}M_S \subseteq N_S$. If $\frac{m}{s} \in N_S$ then $\frac{m}{s} = \frac{n}{t}$ for some $n \in N$ and $t \in S$ so we get $qtm = qsn \in N$ for some $q \in S$ from which we get $m \in N : qt = N$ and if $\frac{r}{s}M_S \subseteq N_S$ then we have $(rM)_S \subseteq N_S$. Now, let $x \in M$. Then $rx \in rM$ and hence for any $s \in S$ we have $\frac{rx}{s} \in (rM)_S \subseteq N_S$ then $\frac{rx}{s} = \frac{n}{t}$ for some $n \in N$ and $t \in S$ thus $qtrx = qsn \in N$ for some $q \in S$ and hence $rx \in N : qt = N$ therefore $rM \subseteq N$. Hence $M$ is a prime submodule of $N$. ■

Definition 3.14. Let $M$ be an $R-$module. If $N$ is a submodule of $M$ and $S$ is a non empty subset of $R$. We define $N : S = \{x \in M : Sx \subseteq N\}$.

As a special case we have $M : S = M$ for each multiplicative closed set $S$ in $R$ and if $N = \{0_M\}$ then we have $\{0_M\} : S = \{x \in M : Sx \subseteq \{0_M\}\} = \{x \in M : Sx = \{0_M\}\} = \{x \in M : sx = 0_M$, for all $s \in S\}$.

Example 3.15. If we take the zero submodule $\{0\}$ of the $Z-$module $Z_6$ and the multiplicative closed set $S = \{-1,1\}$ in $Z$ then $\{0\} : S = \{x \in Z_6 : Sx = \{0\}\} = \{0\}$ and if we take the multiplicative closed set $S = \{2,4\}$ in the $Z_6-$module $Z_6$, then we have $(\{0\} : S) = \{x \in Z_6 : Sx = \{0\}\} = \{0,3\} \neq \{0\}$.

Lemma 3.16. Let $M$ be an $R-$module and $S$ is a multiplicative close set in $R$. If $N'$ is a submodule of $M_S$ the there exists a submodule $N$ of $M$ such that $N' = N_S$. Furthermore, $N$ is a proper submodule of $M$ if and only if $N_S$ is a proper submodule of $M_S$.

Proof. Since $S \neq \phi$, let $s \in S$ be any element and let $N = \{x \in M : \frac{sx}{s} \in N'\}$ (note that $\frac{sx}{s} = \frac{tx}{t}$ for all $s,t \in S$, so there is no any confusion if we take any other element such $t \in S$ to define $N$). We will show $N$ is a submodule of $M$. Clearly $0 = \frac{s0}{s} \in N'$, so that $\phi \neq N \subseteq M$. Let $x, y \in N$ and $r \in R$. Then we have $\frac{s(x-y)}{s} = \frac{sx-sy}{s} = \frac{sx}{s} - \frac{sy}{s} \in N'$, so that $x - y \in N$ and $\frac{s(rx)}{s} = \frac{rs}{s}\frac{sx}{s} \in N'$, so that $rx \in N$ and thus $N$ is a submodule of $M$. To show that $N' = N_S$. Let $x' \in N'$, then $x' = \frac{x}{u}$, for some $x \in M$ and $u \in S$. Now, we have $\frac{sx}{s} = \frac{u(sx)}{us} = \frac{us}{s}\frac{x}{u} = \frac{us}{s}x' \in N'$ and thus $x \in N$, so that $x' = \frac{x}{u} \in N_S$ and hence $N' \subseteq N_S$ and if $x' \in N_S$, then $x' = \frac{n}{v}$ for some $n \in N$ and $v \in S$, then clearly $\frac{sn}{s} \in N'$, so that $x' = \frac{n}{v} = \frac{1}{v}\frac{sn}{s} \in N'$ and thus $N_S \subseteq N'$. Hence $N' = N_S$. Now, let $N \neq M$. If $N_S = M_S$, then for any $x \in M$, we have $\frac{x}{s} \in N_S$, so that $\frac{x}{s} = \frac{a}{l}$, for some $a \in N$ and $l \in S$, thus $\frac{sa}{s} \in N'$, then $\frac{sx}{s} = \frac{s}{s}\frac{sx}{s} = \frac{ss}{s}\frac{x}{s} = \frac{ss}{s}\frac{a}{l} = \frac{s}{l}\frac{sa}{s} \in N'$ and thus $x \in N$, so that $N = M$, which is a contradiction so $N_S \neq M_S$. Now, let $N_S \neq M_S$ and if $N = M$ then clearly $N_S = M_S$ which is a contradiction and thus $N \neq M$. ■

Proposition 3.17. Let $M$ be an $R-$module and $S$ is a multiplicative closed set in $R$ with $(0 : S) = \{0\}$. Then a submodule $N$ of $M$ is essential in $M$ if and only if $N_S$ is essential in $M_S$.

Proof. Let $N$ be essential in $M$ and let $L'$ be any submodule of $M_S$ such that $N_S \cap L' = \{0\}$ then by Lemma 3.16, we have $L' = L_S$ for some submodule $L$ of $M$ and thus $N_S \cap L_S = \{0\}$ this implies that $(N \cap L)_S = \{0\}$. Let $x \in N \cap L$, if $s \in S$ then $\frac{x}{s} \in (N \cap L)_S = 0_S$ and thus $\frac{x}{s} = 0$ which gives that $qx = 0$ for some $q \in S$ and thus $x \in (S : 0) = 0$ so that $x = 0$ and $N \cap L = 0$ and as $N$ is essential we get $L = 0$ and hence $L' = L_S = 0$ and thus $N_S$ is essential in $M_S$.

Conversely, let $N_S$ be essential in $M_S$. If $L$ is any submodule of $M$ such that $N \cap L = 0$ then $L_S$ is a submodule of $M_S$ and $N_S \cap L_S = (N \cap L)_S = 0_S = 0$. As $N_S$ is essential in $M_S$ we get $L_S = 0$. Now let $l \in L$ then for $s \in S$ we have $\frac{l}{s} \in L_S = 0$ and hence $\frac{l}{s} = 0$ so that $pl = 0$ for some $p \in S$ which gives that

$l \in S : 0 = 0$ and thus $l = 0$. Hence $L = 0$ and so that $N$ is essential in $M$. ■

Example 3.18. Consider the $Z_6$ – module $Z_6$. The proper submodules of $Z_6$ are $\{0\}$, $\{0,2,4\}$ and $\{0,3\}$. If we take the multiplicative closed set $S = \{1,5\}$, then it is easy to check that $\{0\} : \{1,5\} = \{0\}$, $\{0,2,4\} : \{1,5\} = \{0,2,4\}$ and $\{0,3\} : \{1,5\} = \{0,3\}$. If we take the multiplicative closed set $S = \{1,3\}$, then $\{0\} : \{1,3\} = \{x \in Z_6 : \{1,3\}x = \{0\}\} = \{0\}$.

Proposition 3.19. Let $M$ be an $R$ – module, $S$ be a multiplicative closed set in $R$ and $N$ is a proper submodule of $M$. If for each proper submodule $K$ of $M$, we have $K : s = K$, for all $s \in S$, then

(1) $N << M$ if and only if $N_S << M_S$.

(2) $N$ is a supplemented submodule of $M$ if and only if $N_S$ is a supplemented submodule of $M_S$.

(3) $M$ is a hollow $R$ – module if and only if $M_S$ is a hollow $R$ – module.

(4) $M$ is a lifting $R$ – module if and only if $M_S$ is a lifting $R$ – module.

Proof. (1) Let $N << M$. Suppose that $L'$ is a submodule of $M_S$ such that $N_S + L' = M_S$, then by Lemma 3.16, $L' = L_S$ for some submodule $L$ of $M$ and so $N_S + L_S = M_S$ which gives that $(N + L)_S = M_S$. If $x \in M$ then for $s \in S$ (since $S \neq \phi$) we have $\frac{x}{s} \in M_S$ and thus $\frac{x}{s} = \frac{n+l}{t}$ for some $n \in N, l \in L$ and $t \in S$, so we get $qtx = qsn + qsl \in N + L$ and thus $x \in (N + L) : p = N + L$. Hence $M \subseteq N + L$. So $N + L = M$. As $N$ is small we get $L = M$ and thus $L' = L_S = M_S$ so that $N_S << M_S$.

Conversely, let $N_S << M_S$. If $L$ is a submodule of $M$ with $N + L = M$ then $N_S + L_S = (N + L)_S = M_S$. As $N_S << M_S$ we get $L_S = M_S$. Now, if $m \in M$, then for any $s \in S$, we have $\frac{m}{s} \in L_S$, that means $\frac{m}{s} = \frac{l}{t}$, for some $l \in L$ and $t \in S$ then for some $q \in S$ we have $qtm = qsl \in L$, thus $m \in L : qt = L$, so that $L = M$ and thus $N << M$.

(2) Let $N$ be a supplemented submodule of $M$, so that $N$ is a supplement of some submodule $L$ of $M$ and thus $N + L = M$ and $N \cap L << N$. Then clearly $N_S + L_S = (N + L)_S = M_S$ and by part (1) we get $N_S \cap L_S = (N \cap L)_S << N_S$, so that $N_S$ is a supplement of $L_S$ in $M_S$ and thus $N_S$ is a supplemented submodule of $M_S$.

(3) Let $M$ be a hollow $R$ – module and $N'$ be any proper submodule of $M_S$, then $N' = N_S$, for some submodule $N$ of $M$. As $N_S$ is proper in $M_S$ we get $N$ is proper in $M$ so that $N << M$ and thus by part (1) we get $N' = N_S << M_S$. Hence $M_S$ is a hollow $R$ – module.

(4) Let $M$ be a lifting $R$ – module. To show $M_S$ is a lifting $R$ – module, let $N'$ be any submodule of $M_S$, so that $N' = N_S$ for some submodule $N$ of $M$. As $M$ is lifting there exist submodules $K, L$ of $M$ such that $M = K \oplus L$ with $K \leq N$ and $N \cap L << L$. Clearly, $M = K + L$, so that $M_S = (K + L)_S = K_S + L_S$ and if $\frac{x}{s} \in K_S \cap L_S$, where $x \in M, s \in S$ then $\frac{x}{s} = \frac{k}{t} = \frac{l}{u}$ for some $k \in K, l \in L$ and $t, u \in S$. Hence there exist $p, q \in S$ such that $ptx = psk \in K$ and $qux = qsl \in L$, thus $pqutx \in K$ and $pqutx \in L$, so that $pqutx \in K \cap L = \{0\}$ and then $pqutx = 0$, $\frac{x}{s} = \frac{pqutx}{pquts} = 0$, so that $K_S \cap L_S = \{0\}$. Hence $M_S = K_S \oplus L_S$ and also we have $K_S \leq N_S$ and $N_S \cap L_S = (N \cap L)_S << L_S$. Hence $M_S$ is a lifting $R$ – module.

Conversely, suppose that $M_S$ is a lifting $R$ – module. Let $N$ be any submodule of $M$ then $N_S$ is a submodule of $M_S$. Hence there exist submodules $K_S, L_S$ of $M_S$, where $K, L$ are submodules of $M$, such that $M_S = K_S \oplus L_S$ with $K_S \leq N_S$ and $N_S \cap L_S << L_S$. We proceed as in the first part and we get $M = K \oplus L$ with $K \leq N$ and $N \cap L << L$. Hence $M$ is a lifting $R$ – module. ■

Example 3.20. In the $Z$ – module $Z_6$, the only proper submodules of $Z_6$ are $\{0\}, \{0,2,4\}$ and $\{0,3\}$. Consider the multiplicative closed set $S = \{1,5\}$ in $Z_6$. We have $S(\{0\}) = \{m \in Z_6 : mx \in \{0\}$ for some $x \notin \{0\}\} = \{m \in Z_6 : mx = 0$, for some $x \neq 0\} = \{0,2,3,4\}$ (which is not an ideal of $Z_6$ and so that $\{0\}$ is not a primal ideal of $Z_6$) and clearly we have $S(\{0\}) \cap S = \phi$. Also we have $S(\{0,2,4\}) = \{m \in Z_6 : mx \in \{0,2,4\}$ for some $x \notin \{0,2,4\}\} = \{m \in Z_6 : mx \in \{0,2,4\}$, for $x \in \{1,3,5\}\} = \{0,2,4\}$ and clearly $S(\{0,2,4\}) \cap S = \phi$. Finally, we have $S(\{0,3\}) = \{m \in Z_6 : mx \in \{0,3\}$ for some $x \notin \{0,3\}\} = \{m \in Z_6 : mx \in \{0,3\}$, for $x \in \{1,2,4,5\}\} = \{0,3\}$ and clearly $S(\{0,3\}) \cap S = \phi$, so that $S(N) \cap S = \phi$, for every proper submodule of the $Z_6$ – module $Z_6$.

Lemma 3.21. Let $M$ be an $R$ – module and $S$ is a multiplicative closed set in $R$ such that $S(K) \cap S = \phi$, for every submodule $K$ of $M$, a submodule $N$ is a maximal submodule of $M$ if and only if $N_S$ is a maximal submodule of $M_S$.

Proof. Suppose that $N$ is a maximal submodule of $M$. To show $N_S$ is a maximal submodule of $M_S$. As $N$ is a proper submodule of $M$, we have $N_S$ is a proper submodule of $M_S$. Let $L'$ be any submodule of $M_S$ such that $N_S \subseteq L' \subseteq M_S$, then by part (1), we have $L' = L_S$ for some submodule $L$ of $M$ and so that we have $N_S \subseteq L_S \subseteq M_S$. To show $N \subseteq L$. Let $x \in N$, then $\frac{x}{s} \in L_S$, so that $\frac{x}{s} = \frac{l}{t}$, for $l \in L$ and $t \in S$, then $ptx = psl \in L$, for some $p \in S$. If $x \notin L$, then $pt \in S(L)$ and as $S(L) \cap S = \phi$, so that $pt \notin S$ which is a contradiction, and thus $x \in L$. Hence $N \subseteq L \subseteq M$ and as $N$ is maximal so that $N = L$ or $L = M$. This implies $N_S = L_S$ or $L_S = M_S$, so that $N_S$ is a maximal submodule of $M_S$. The converse part can be done similarly. ■

Lemma 3.22. Let $M$ be an $R$ – module and $S$ is a multiplicative closed set in $R$ such that $S(K) \cap S = \phi$, for every submodule $K$ of $M$, then

(1) If $N$ is any submodule of $M$, then $\text{Rad}\, N_S = (\text{Rad}\, N)_S$.

(2) If $N$ and $L$ are submodules of $M$, then $N = L$ if and only if $N_S = L_S$.

(3) $P(M_S) = (P(M))_S$.

Proof. (1) Let $a' \in \text{Rad}\, N_S$. Then $a' = \frac{a}{s}$, for $a \in N$ and $s \in S$. Let $L$ be any maximal submodule of $N$ then by Lemma 3.21, $L_S$ is a maximal submodule of $N_S$ and thus $a' = \frac{a}{s} \in L_S$, then $\frac{a}{s} = \frac{l}{t}$, for some $l \in L$ and $t \in S$, thus $pta = psl \in L$. If $a \notin L$ then $pt \in S(L)$ and as $S(L) \cap S = \phi$, we get $pt \notin S$ which is a contradiction so that $a \in L$ and so $a \in \text{Rad}\, N$, so that $a' = \frac{a}{s} \in (\text{Rad}\, N)_S$. Hence $\text{Rad}\, N_S \subseteq (\text{Rad}\, N)_S$.

Conversely, suppose that $x' \in (\text{Rad}\, N)_S$, so that $x' = \frac{x}{u}$, for some $x \in \text{Rad}\, N$ and $u \in S$. Let $L'$ be any maximal submodule of $N_S$, then by Lemma 3.16, we have $L' = L_S$, for some submodule $L$ of $N$ and as $L_S$ is maximal by Lemma 3.21, we get $L$ is a maximal submodule of $N$, so that $x \in L$ and so $x' = \frac{x}{u} \in L_S = L'$ and hence $x' \in \text{Rad}\, N_S$, so that $(\text{Rad}\, N)_S \subseteq \text{Rad}\, N_S$. Hence $\text{Rad}\, N_S = (\text{Rad}\, N)_S$.

(2) Let $N = L$, then clearly $N_S = L_S$. Conversely, let $N_S = L_S$ to show that $N = L$. Let $n \in N$, then for an $s \in S$ ($S \neq \phi$), we have $\frac{n}{s} \in L_S$, so that $\frac{n}{s} = \frac{l}{t}$, for some $l \in L$ and $t \in S$, then $ptn = psl \in L$, for $p \in S$ and if $n \notin L$, then $pt \in S(L)$ and as $S(L) \cap S = \phi$ we must have $pt \notin S$ which is a contradiction, and so $n \in L$, that gives $N \subseteq L$ and by the same way we get $L \subseteq N$. Hence $N = L$.

(3) Let $x' \in P(M_S) = \sum_{\alpha \in \Delta} N'_\alpha$, where $N'_\alpha$ is a submodule of $M_S$ with $\text{Rad}\, N'_\alpha = N'_\alpha$, for each $\alpha \in \Delta$. By Lemma 3.16, for each $\alpha \in \Delta$ we have $N'_\alpha = (N_\alpha)_S$, where $N_\alpha$ is a submodule of $M$, and thus $\text{Rad}\,((N_\alpha)_S) = (N_\alpha)_S$, then by what we have shown in part (2), we get $(\text{Rad}\, N_\alpha)_S = \text{Rad}((N_\alpha)_S) = (N_\alpha)_S$, and by part(2), we get $\text{Rad}\, N_\alpha = N_\alpha$, for each $\alpha \in \Delta$. As $x' \in \sum_{\alpha \in \Delta} N'_\alpha$, there exists a positive integer $m$ with a finite subset $J = \{\alpha_i\}_{i=1}^{m}$ of $\Delta$ such that $x' = \sum_{i=1}^{m} n'_{\alpha_i}$, where $n'_{\alpha_i} \in N'_{\alpha_i}$ and $N'_{\alpha_i} = (N_{\alpha_i})_S$, so that $\text{Rad}\,((N_{\alpha_i})_S) = (N_{\alpha_i})_S$, $(\text{Rad}\, N_{\alpha_i})_S = \text{Rad}((N_{\alpha_i})_S) = (N_{\alpha_i})_S$ and $\text{Rad}\, N_{\alpha_i} = N_{\alpha_i}$. Now, for each $i$ $(1 \le i \le m)$, we have $n'_{\alpha_i} \in (N_{\alpha_i})_S$ and thus for each $i$ $(1 \le i \le m)$, we have $n'_{\alpha_i} = \frac{n_{\alpha_i}}{s_i}$, where $n_{\alpha_i} \in N_{\alpha_i}$ and $s_i \in S_i$, then we have $x' = \sum_{i=1}^{m} n'_{\alpha_i} = \sum_{i=1}^{m} \frac{n_{\alpha_i}}{s_i} =$

$\frac{\sum_{i=1}^{m} s_1 s_2 .. s_{i-1} s_{i+1} .. s_m n_{\alpha_i}}{s_1 s_2 .. s_m} \in (\sum_{i=1}^{m} N_{\alpha_i})_S \subseteq (P(M))_S$, so that $P(M_S) \subseteq (P(M))_S$. In a similar argument we can prove that $(P(M))_S \subseteq P(M_S)$. Hence $P(M_S) = (P(M))_S$. ■

Remark 3.23. If $S$ is a multiplicative closed set in $R$, then for any proper submodule $N$ with $S(N) \cap S = \phi$, then for any $p \in S$, we have $N : p \subseteq N$. Let $x \in N : p$ and then $px \in N$. If $x \notin N$ then we have $p \in S(N)$, so that $p \notin S$ which is a contradiction, so that $x \in N$ and so that $N : p \subseteq N$. Hence $N : p = N$.

Proposition 3.24. Let $M$ be an $R-$module and $S$ is a multiplicative closed set in $R$ such that $S(K) \cap S = \phi$, for every proper submodule $K$ of $M$, then

(1) $M$ is coatomic if and only if $M_S$ is coatomic.

(2) $M$ is reduced if and only if $M_S$ is reduced.

(3) $M$ is a local $R-$module if and only if $M_S$ is a local $R_S$ $-$module.

Proof. (1) Let $M$ be a coatomic module and $N'$ be any proper submodule of $M_S$, then by Lemma 3.16, we have $N' = N_S$ for some submodule $N$ of $M$ and as $N_S$ is a proper submodule of $M_S$, by Lemma 3.16, we have $N$ is a proper submodule of $M$. Then $N \subseteq L$, for some maximal submodule $L$ of $M$. Then by Lemma 3.21, we have $L_S$ is a maximal submodule of $M_S$ with $N' = N_S \subseteq L_S$ and thus $M_S$ is coatomic. Conversely, suppose that $M_S$ is a coatomic $R-$module and let $N$ be a proper submodule of $M$, then by Lemma 3.16, $N_S$ is a proper submodule of $M_S$, so that $N_S \subseteq L'$, where $L'$ is a maximal submodule of $M_S$ and then by Lemma 3.16, we have $L' = L_S$ for some maximal submodule $L$ of $M$ with $N_S \subseteq L_S$ and by what we have shown in part (2) of Lemma 3.22, we get $N \subseteq L$, so that $M$ is a coatomic $R-$module.

(2) Let $M$ be reduced, so that $P(M) = 0$ and by part (3) of Lemma 3.22, we have $P(M_S) = (P(M))_S = 0$. Conversely, suppose that $P(M_S) = 0$ and to show $P(M) = 0$. Let $x \in P(M)$, then for any $s \in S$ we have $\frac{x}{s} \in (P(M))_S = P(M_S) = 0$ and so that $qx = 0 \in \{0\}$, for some $q \in S$. If $x \neq 0$, so that $x \notin \{0\}$ and thus $q \in S(\{0\})$ and as $S(\{0\}) \cap S = \phi$, we must have $q \notin S$ which is a contradiction so that $x = 0$. Hence $P(M) = 0$.

(3) Suppose that $M$ is a local $R-$module. To show $M_S$ is a local $R_S$ $-$module. By Lemma 3.22, we have $\mathrm{Rad}\, M_S = (\mathrm{Rad}\, M)_S$. As $M$ is a local, we have $\mathrm{Rad}\, M$ is a maximal submodule of $M$ and $\mathrm{Rad}\, M << M$. Then by Lemma 3.21, we have $(\mathrm{Rad}\, M)_S$ is a maximal submodule of $M_S$ and by the above remark $(\mathrm{Rad}\, M) : s = \mathrm{Rad}\, M$, for all $s \in S$, so by Proposition 3.19, we have $(\mathrm{Rad}\, M)_S << M_S$, that is, $\mathrm{Rad}\, M_S$ is a maximal submodule of $M_S$ and $\mathrm{Rad}\, M_S << M_S$ and thus $M_S$ is a local $R_S$ $-$module. The proof of the converse part can be done similarly. ■

Remark 3.25. If $P$ is a prime ideal of $R$, then $R-P$ is a multiplicative closed set in $R$, then for any proper submodule $N$ with $S(N) \subseteq P$ we have

(1) $S(N) \cap (R-P) = \phi$.

(2) If $p \notin P$, then to show $N : p \subseteq N$, so let $x \in N : p$ and then $px \in N$. If $x \notin N$ then we have $p \in S(N) \subseteq P$, which is a contradiction so that $x \in N$ and so that $N : p \subseteq N$. Hence $N : p = N$, so that we can give the following corollary.

Corollary 3.26. Let $M$ be an $R-$module and $P$ is a prime ideal of $R$ such that $S(K) \subseteq P$, for each proper submodule $K$ of $M$, then:

(1) $M$ is coatomic if and only if $M_P$ is coatomic.

(2) $M$ is reduced if and only if $M_P$ is reduced.

(3) $M$ is a hollow $R-$module if and only if $M_P$ is a hollow $R-$module.

(4) $M$ is a lifting $R-$module if and only if $M_P$ is a lifting $R-$module.

(5) $M$ is a local $R-$module if and only if $M_S$ is a local $R_S$ $-$module.

(6) A submodule $N$ of $M$ is a maximal submodule of $M$ if and only if $N_P$ is a maximal submodule of $M_P$.

(7) If $N$ is a proper submodule of $M$ then $S(N)$ is an ideal of $R$ if and only if $S(N_P)$ is an ideal of $R_P$.

(8) If $N$ is a proper submodule of $M$ then $N$ is primal if and only if $N_P$ is primal.

(9) If $N$ is a proper submodule of $M$ then $N$ is essential in $M$ if and only if $N_P$ is essential in $M_P$.

(10) If $N$ is a proper submodule of $M$ then $N << M$ if and only if $N_P << M_P$.

(11) If $N$ is a proper submodule of $M$ then $N$ is a supplemented submodule of $M$ if and only if $N_P$ is a supplemented submodule of $M_P$.

(12) If $N$ is a proper submodule of $M$ then Rad $N_P = (\text{Rad } N)_P$.

Remark 3.27. If $R$ is a local ring with the unique maximal ideal $P$ and if $N$ is any primal submodule of $M$, then $S(N)$ is an ideal of $R$ and thus $S(N) \subseteq P$.

Corollary 3.28. Let $R$ be a local ring with the unique maximal ideal $P$ and $M$ is an $R-$module. If $N$ is a primal submodule of $M$, then:

(1) $M$ is coatomic if and only if $M_P$ is coatomic.

(2) $M$ is reduced if and only if $M_P$ is reduced.

(3) $M$ is a hollow $R-$module if and only if $M_P$ is a hollow $R-$module.

(4) $M$ is a lifting $R-$module if and only if $M_P$ is a lifting $R-$module.

(5) $M$ is a local $R-$module if and only if $M_P$ is a local $R_P$ $-$module.

(6) $N$ is a maximal submodule of $M$ if and only if $N_P$ is a maximal submodule of $M_P$.

(7) $N$ is an essential submodule of $M$ if and only if $N_P$ is essential in $M_P$.

(8) $N << M$ if and only if $N_P << M_P$.

(9) $N$ is a supplemented submodule of $M$ if and only if $N_P$ is a supplemented submodule of $M_P$.

(10) Rad $N_P = (\text{Rad } N)_P$.

Remark 3.29. let $R$ be a commutative ring with identity and $S$ is the set of all units of $R$, then clearly $0 \notin S$ and for $u, v \in S$ we have $uv \in S$, so that $S$ is a multiplicative closed set in $R$, also we have $1 \in S$. Now, let $M$ be an $R-$module and $N$ is a primal submodule of $M$. If $S(N) \cap S \neq \phi$ then there exists $s \in S(N)$ and $s \in S$, so that $s$ is a unit and as $S(N)$ is an ideal of $R$, we get $S(N) = R$ which is a contradiction. Hence $S(N) \cap S = \phi$. Hence we can give the following corollary.

Corollary 3.30. Let $M$ be an $R-$module and $S$ is the set of all units of $R$. If $N$ is a primal submodule of $M$, then

(1) $M$ is coatomic if and only if $M_S$ is coatomic.

(2) $M$ is reduced if and only if $M_S$ is reduced.

(3) $M$ is a hollow $R-$module if and only if $M_S$ is a hollow $R-$module.

(4) $M$ is a lifting $R-$module if and only if $M_S$ is a lifting $R-$module.

(5) $M$ is a local $R-$module if and only if $M_S$ is a local $R_S$ $-$module.

(6) $N$ is a maximal submodule of $M$ if and only if $N_S$ is a maximal submodule of $M_S$.

(7) $N$ is essential in $M$ if and only if $N_S$ is essential in $M_S$.

(8) $N << M$ if and only if $N_S << M_S$.

(9) $N$ is a supplemented submodule of $M$ if and only if $N_S$ is a supplemented submodule of $M_S$.

(10) Rad $N_S = (\text{Rad } N)_S$.

References.

## بعض الشروط والتي تجعل انواعا محددة من المقاسات تملك خاصية التوضيع

اسماعيل عقراوي [1] , عادل قادر جبار [2] و رضا سزيده [3]

[1] قسم الؤياضيات - فاكلتي العلوم – جامعة سوران

[2] قسم الرياضيات – سكول العلوم – فاكلتي العلوم و تربية العلوم – جامعة السليمانية

[3] قسم الرياضيات - فاكلتي العلوم – جامعة ورمى

### الخلاصة

نعني بخاصية التوضيع الخاصية التي تبقى تحت تأثير عملية توضيع المقاسات عند المجموعات المغلقة ضربيا. الهدف من هذا البحث هو ايجاد الشروط التي تمكننامن نقل خاصية معينة لمقاس معلوم الى توضيعه عند المجموعات المغلقة ضربيا و بالعكس وهذا يعني تحديد الشروط التي تجعل من خاصية معلومة لمقاس لتكون خاصية التوضيع ولذلك نعطي شروطا عديدة والتي عند توفرها تجعل انواعا محددة من المقاسات تملك خاصية التوضيع.

## هەندێ مەرج کەوا لە چەند جۆرێکی دیاری کراو لە پێوەرەکان دەکات کە تایبە تمەندی خوجێیتیان هەبێت

اسماعيل عقراوي [1] , عادل قادر جبار [2] و رضا سزيده [3]

[1] بەشی ما تماتیك-فاكەلتی زانست – زانكۆی سۆران

[2] بەشی ما تماتیك- سكولی زانست - فاكەلتی زانست و پەروەردە زانستەكان– زانكۆی سلێمانی

[1] بەشی ما تماتیك-فاكەلتی زانست – زانكۆی ورمێ

### كورتە

مەبەستمان لە تايبە تمەندی خۆجێیتی ئەو تايبە تمەندييە كە دەمێنێتەوە لە ژێر كاريگەری بەخۆجێ كردنی پێوەرەكان لە كۆمەڵە بەلێكدان داخراوەكان. ئامانج لەم توێژينەوەيە دۆزينەوەی ئەو مەرجانەيە كە بە هۆيانەوە دەتوانين تايبە تمەندێكی دياری كراوی پێوەرێكی دراو بگوێزينەوە بۆ خۆجێتيەكەی ە كۆمەڵە بەلێكدان داخراوەكان و بە پێچەوانەشەوە ئەمەش واتای ئەوەيە كە دەتوانين ئەو مەرجانە دياری بكەين كە وا لە تايبە تمەندێكی زانراو ی پێوەرێك دەكات تايبە تمەندی خۆجێيتی هەبێت و هەروەها چەند مەرجێك دەدەين كە بە هۆيانەوە چەند پێوەرێك دەتوانێت ببێت بە خاوەنی تايبە تمەند خۆجێيتی.